\documentclass[11pt,a4paper]{amsart}   
\usepackage[latin1]{inputenc}
\usepackage[T1]{fontenc}               
\usepackage[francais]{babel}
\usepackage{amssymb}
\usepackage{mathrsfs}
\newtheorem{theo}{Th\'eor\`eme}[section]
\newtheorem{lemme}[theo]{Lemme}
\newtheorem{prop}[theo]{Proposition}
\newtheorem{cor}[theo]{Corollaire}
\newtheorem{rem}[theo]{Remarque }
\newenvironment{remarque}{\begin{rem}\em}{\end{rem}}
\newtheorem{rems}[theo]{Remarques }
\newenvironment{remarques}{\begin{rems}\em}{\end{rems}}
\newtheorem{defi}[theo]{D\'efinition } 

\newenvironment{dem}{\noindent{\it D\'emonstration}.}
{\unskip\hfill\null\nobreak\hfill\carre\vskip1em\par}
\newcommand{\carre}{\rule{1ex}{1ex}}
\newcommand{\rond}{\raisebox{.3mm}{\scriptsize$\circ$}}

\def\g{\mathfrak}
\def\R{\mathbb R}
\def\C{\mathbb C}
\def\H{\mathscr H}
\def\A{\mathscr A}
\def\E{\mathscr E}

\def\Sc{\mathscr S}
\def\D{\mathscr D}
\def\P{\mathscr P}
\def\Ad{\mathop{\mbox{\rm Ad}}\nolimits}
\def\ad{\mathop{\mbox{\rm ad}}\nolimits}
\def\tr{\mathop{\mbox{\rm tr}}\nolimits}
\title{Sur les distributions covariantes dans les algèbres de Lie
  réductives}
\author{Abderrazak Bouaziz}

\begin{document}
\maketitle

\begin{abstract}Nous étudions le problème de la décomposition d'une
  distribution covariante sur une algèbre de Lie réductive en somme
  finie de produits d'une distribution invariante par un polyn\^ome covariant.

\end{abstract}

\section{Introduction}

Soit $G$ un groupe de Lie agissant linéairement dans un espace
vectoriel réel $E$ et dans un espace vectoriel réel ou complexe
$V$. Une distribution covariante sur $E$ à valeurs dans $V$ est une application
linéaire continue $G$-équivariante de $C_c^{\infty}(E)$ dans $V$. Par
exemple si  $P$ est un
 polyn\^ome $G$-équivariant de $E$ dans $V$, on associe à chaque
 distribution $G$-invariante $\theta $ sur $E$ une distribution
 $G$-équivariante $\theta P$ définie de la façon suivante.  On écrit
 $P=\sum_ip_i\otimes v_i$, où $(v_1,\ldots,v_n)$ est une base de $V$, et, pour
 $f\in C_c^{\infty}(E)$, on pose
$$\langle \theta P,f\rangle=\sum_i\langle \theta ,p_if\rangle v_i.$$
Cette définition ne dépend pas du choix de la base de $V$.

Si $G$ est compact et l'action de $G$ dans $E$ est colibre, A.I. Oksak
\cite{ok}  a montré que toute distribution covariante $T$ se factorise en
somme finie $T=\theta_1P_1+\cdots+\theta_kP_k$, où $P_1,\ldots, P_k$
est une base de l'espace des polyn\^omes covariants comme module sur
l'anneau des polyn\^omes invariants sur $E$ ; ce résultat a été étendu
récemment par A. Saidi \cite{Sa} au cas général de l'action d'un
groupe compact (sans l'hypothèse colibre).

Pour $G=SL(2,\R)$, P. Lavaud \cite{lav} a décrit les distributions
covariantes sur $sl(2,\R)$ à valeurs dans $sl(2,\R)$ dont le  support
est inclus dans
le c\^one nilpotent. Il a montré en particulier, mais avec une
formulation différente, que de telles  distributions se factorisent sous la
forme $\theta P$, où $P$ est un polyn\^ome covariant fixe et $\theta $
parcourt l'ensemble des distributions $G$-invariantes sur $\g g$ à
support inclus dans le c\^one nilpotent. 

Cet exemple constitue la motivation principale de ce travail. Nous nous proposons d'entreprendre l'étude des
distributions covariantes sur une algèbre de Lie réductive $\g g$ à
valeurs dans l'espace d'une représentation de dimension finie $V$ du
groupe adjoint $G$ de $\g g$. Dans cette situation, on sait, d'après Kostant \cite{Kos}, qu'il
existe un nombre fini de polyn\^omes covariants de $\g g$ dans $V$
tels que  tout polyn\^ome covariant $P$ de $\g g$ dans $V$ se
factorise de façon unique en somme $P=Q_1P_1+\cdots+Q_rP_r$, où les
$Q_i$  sont des polyn\^omes invariants. On se pose alors la question
suivante.

{\it Est-ce que toute distribution covariante $T$ dans  $\g g$  à
  valeurs dans $V$ se factorise sous la forme
  $T=\theta_1P_1+\cdots+\theta_rP_r$, où les $\theta_i$ sont des
  distributions invariantes ?}

Dans cet article, nous montrons une forme faible
de cette décomposition (où les $\theta_i$ ne sont pas forcément invariantes) pour
les distributions covariantes sur $\g g$ à valeurs dans $\g g$ ; notre démonstration  s'inspire beaucoup
des méthodes développées dans \cite{Dix}. Nous montrons aussi le
m\^eme type de décomposition pour les fonctions lisses et analytiques,
et nous étendons au cas lisse la description  de Dixmier des champs de
vecteurs tangents aux orbites, répondant ainsi
à une question de Miranda et Zung  \cite{mir}.

\section{Notations}\label{notations}

Soit $\g g$ une algèbre de Lie réductive réelle et soit $G$ son groupe adjoint.
On note $\ell$ le rang de $\g g$. Si $x\in \g g$, soit $\g g^x$
(resp. $G^x$) son centralisateur dans $\g g$ (resp. $G$). On dit que
$x\in \g g$ est régulier si $\dim \g g^x=\ell$. On note $\g r$ l'ensemble des
éléments réguliers de $\g g$.

On fixe une forme bilinéaire symétrique non dégénérée $G$-invariante
$\kappa$ sur $\g g$.

On note $\g g_{\C}$ l'algèbre de Lie complexifiée de $\g g$,
$G_{\C}$ son groupe adjoint  et $\g r_{\C}$ l'ensemble de ses 
éléments réguliers.

On fixe  une représentation de dimension finie  $(\pi, V)$ de $G_{\C}$
($V$ espace vectoriel sur $\C$). On
note encore $\pi$ la différentielle de $\pi$. Les poids de $V$ sont
dans le réseau des racines ; en particulier $0$ est un poids de $V$,
on note $r$ sa multiplicité. Si $H$ est un sous-groupe de $G_{\C}$, on
note $V^{H}$ l'ensemble des points fixes de $H$ dans $V$.

On note $\mathcal P(\g
 g_{\C},V)^{G_{\C}}$ l'espace des polynômes $G_{\C}$-équivariants sur 
 $\g g_{\C}$  à valeurs dans $V$, c'est-à-dire les polyn\^omes $P:\g g_{\C}\rightarrow V$
 vérifiant
$$P(\Ad(g)\cdot x)=\pi(g)\cdot P(x),\quad \mbox{\rm pour tous  } g\in
 G_{\C}, x\in \g g_{\C},$$
et on note $\mathcal P(\g g_{\C})^{G_{\C}}$ l'algèbre des polyn\^omes
 $G_{\C}$-invariants sur $\g g_{\C}$. D'après Kostant (\cite{Kos}, théorème
 9) :

 K1. $P(\g g_{\C}, V)^{G_{\C}}$ est un $\mathcal P(\g
 g_{\C})^{G_{\C}}$-module libre de rang $r$,

 K2.   si $P_1,\ldots,P_r$ est une base de $\mathcal P(\g
   g_{\C},V)^{G_{\C}}$ formée d'éléments homogènes,  alors pour tout $x\in \g
   r_{\C}$, les vecteurs $P_1(x),\ldots,P_r(x)$ forment une base de
   $V^{G_{\C}^x}$.

Dans \cite{Kos} l'algèbre de Lie est supposée semi-simple, mais
l'extension de ces résultats au cas réductif est évidente.

\begin{remarque}\label{formereelle}{Si $V$ admet une forme réelle $V_{\R}$ stable par
   $\pi(G)$, on peut prendre les $P_i$ de sorte que $P_i(\g g)\subset
   V_{\R}$ (voir appendice).}
\end{remarque}

Si $U$ est un ouvert d'un espace vectoriel réel de dimension finie et
si $F$ est
un espace de Fréchet sur $\R$ ou sur $\C$, on note $\E(U,F)$   l'espace des
fonctions de classe $C^{\infty}$ de $U$ dans $F$, et
$\E_c(U,F)$ le sous espace de $\E(U,F)$ formé des
fonctions à support compact. Sauf mention explicite du contraire, la
notation $\E(U)$ désignera $\E(U,\C)$. Si $K$ est une
partie compacte de $U$, on note $\E_K(U,F)$ le sous-espace de
$\E_c(U,F)$ formé des fonctions à support inclus dans
$K$.  Tous ces espaces seront munis de leurs topologies naturelles.

\section{Fonctions covariantes}\label{secun}

D'après les propriétés K1 et K2 ci-dessus, tout polynôme $P\in \mathcal P(\g
 g_{\C},V)^{G_{\C}} $ s'écrit de façon unique
 $P~=~R_1P_1+\cdots+R_rP_r$, avec $R_1,\ldots,R_r\in \mathcal P(\g
 g_{\C})^{G_{\C}}$. On se propose de montrer des décompositions
 analogues dans différents espaces de fonctions.

\subsection{Fonctions holomorphes covariantes}

Si $U$ est un ouvert de $\g g_{\C}$, on  note  $\H(U,V)$ l'espace des
fonctions holomrphes de $U$ dans $V$ ; si $V=\C$, on le note
simplement $\H(U)$. Si de plus $U$ est
$G_{\C}$-invariant, on note $\H(U,V)^{G_{\C}}$ le sous-espace de
$\H(U,V)$ formé des fonctions $G_{\C}$-covariantes.

\begin{prop}\label{hol}{Soient $U$ un ouvert de $\g g_{\C}$ et $f\in \H(U,V)$.  Alors les   conditions suivantes sont équivalentes :

{\em(i)} pour tout $x\in U$, on a $f(x)\in V^{G_{\C}^x}$ ;

{\em(ii)} pour tout $x\in U\cap \g r_{\C}$, on a $f(x)\in V^{G_{\C}^x}$ ;

{\em(iii)} il existe des applications $f_1,\ldots,f_r\in \H(U)$,
uniquement déterminées, telles
que $f=f_1P_1+\cdots+f_rP_r$.}
\end{prop}
\begin{dem} Il est clair que (iii)$\Rightarrow$(i)$\Rightarrow$(ii). Supposons que $f$
  vérifie (ii). Si $\g g_{\C}$ est abélienne, les $P_i$ sont des
  constantes et tout est trivial. On  suppose donc que  $\g g_{\C}$ n'est pas abélienne. D'après K2, il  existe  $r$ fonctions
  $f_1,\ldots,f_r$ à valeurs complexes sur $U\cap \g r_{\C}$, uniquement déterminées, telles que
  \begin{equation}\label{un}
f(x)=f_1(x)P_1(x)+\cdots+f_r(x)P_r(x)\quad\mbox{\rm pour tout  } x\in U\cap \g
  r_{\C}.
\end{equation}
 Il est clair que les $f_i$ sont holomorphes. D'autre part,
  la codimension de $U\setminus \g r_{\C}$ dans $U$ est $3$
  (\cite{Ve}, th. 4.12). Donc, d'après le théorème de Hartogs, les
  $f_i$ se prolongent en des fonctions holomorphes sur $U$, qu'on note
   de la m\^eme façon. L'identité (\ref{un}) est donc vraie pour
  tout $x\in U$ par densité de $ U\cap \g
  r_{\C}$ dans $U$.
\end{dem}

\begin{cor}\label{holinv}{Soient $U$ un ouvert $G_{\C}$-invariant de $\g g_{\C}$ et
  $f\in \H(U,V)^{G_{\C}}$. Alors il existe $r$ fonctions $f_1,\ldots,f_r\in
  \H(U)^{G_{\C}}$, uniquement déterminées, telles que $f=f_1P_1+\cdots f_rP_r$.}
\end{cor}
\begin{dem}
La covriance implique que $f(x)\in V^{G_{\C}^x}$ pour tout $x\in
U$. Le théorème implique alors l'existence des $f_i\in \H(U)$ ; l'invariance
des $f_i$
découle de leur unicité.
\end{dem}

\begin{remarque}{Ce corollaire est un cas particulier d'un résultat de
    G. Schwarz (\cite{Schz}, proposition 6.8) obtenu par une autre
    méthode.}
\end{remarque}

\subsection{Fonctions analytiques réelles  covariantes}\label{troisdeux} 

Si $U$ est un ouvert de $\g g$, on note $\A(U,V)$ l'espace des
fonctions analytiques de $U$ dans $V$. On note
$\A(U)$ pour $\A(U,\C)$.

On fixe une base $q_1,\ldots,q_{\ell}$ de $P(\g g_{\C},\g g_{\C})^{G_{\C}}$, en
tant que $\mathcal P(\g g_{\C})^{G_{\C}}$-module, formée d'éléments
homogènes et telle que $q_i(\g g)\subset \g g$ pour tout
$i=1,\ldots,\ell$. L'existence d'une telle base découle de  la remarque
\ref{formereelle} ; on peut aussi la voir d'une autre façon : il est
bien connu que $\P(\g g_{\C})^{G_{\C}}$ admet un système générateur
$p_1,\ldots,p_{\ell}$ formé de polyn\^omes homogènes réels sur $\g g$,
alors, si on identifie $\g g$ à $\g g^*$ à l'aide de  $\kappa$, les  polyn\^omes
$dp_1,\ldots, dp_{\ell}$ ont les propriétés requises.

\begin{prop}\label{ana}{Soient $U$ un ouvert de $\g g$ et $f\in \A(U,V)$.  Alors les    conditions suivantes sont équivalentes :

{\em(i)} pour tout $x\in U$, on a $f(x)\in V^{G^x}$ ;

{\em(ii)} pour tout $x\in U\cap \g r$, on a $f(x)\in V^{G^x}$ ;

{\em(iii)} il existe des applications $f_1,\ldots,f_r\in \A(U)$,
uniquement déterminées, telles
que $f=f_1P_1+\cdots+f_rP_r$.}
\end{prop}
\begin{dem} Il s'agit de montrer (ii)$\Rightarrow$(iii).  Soit $\tilde{f}$ un prolongement holomorphe de $f$ à un
  voisinage ouvert  $\tilde{U}$ de $U$ dans $\g g_{\C}$. On peut
  supposer que chaque composante connexe de $\tilde{U}$ rencontre
  $U$. Il suffit alors de
  montrer que $\tilde{f}$ vérifie la condition (ii) de la proposition
  \ref{hol}, car, par restriction à $U$, la décomposition de
  $\tilde{f}$ fournit la décomposition de $f$. 

Soit $x\in \tilde{U}\cap \g r_{\C}$. D'après (\cite{Kos}, lemme 5),
$G_{\C}^x$ est connexe, donc $V^{G_{\C}^x}=V^{\g g_{\C}^x}$. Comme les
vecteurs $q_1(x),\ldots,q_{\ell}(x)$ engendrent $\g g_{\C}^x$ (sur
$\C$), il suffit de montrer que $\pi(q_i(x))\cdot \tilde{f}(x)=0$ pour tout
$i=1,\ldots,\ell$ et pour tout $x\in \tilde{U}\cap \g r_{\C}$. 

La fonction $ g_i:x\mapsto\pi(q_i(x))\cdot \tilde{f}(x)$ est holomorphe sur
$\tilde{U}$. Pour tout $x\in \g g$, on a $V^{G^x}\subset V^{\g
  g^x}$, donc d'après (ii) $g_i$  est nulle sur $U\cap \g r$ et donc
nulle sur $U$. Il s'ensuit qu'elle  est identiquement nulle,
car $U$ rencontre toutes les composantes connexes de $\tilde{U}$. 
\end{dem}
\begin{cor}\label{corana}{Soient $U$ un ouvert $G$-invariant de $\g g$ et
  $f\in \A(U,V)^{G}$. Alors il existe $r$ fonctions $f_1,\ldots,f_r\in
  \A(U)^{G}$, uniquement déterminées, telles que $f=f_1P_1+\cdots f_rP_r$.}
\end{cor}
\begin{dem} Analogue à celle du corollaire \ref{holinv}
\end{dem}

\subsection{Fonctions lisses  covariantes}

\begin{prop}\label{lisse}{Soient $U$ un ouvert de $\g g$ et $f\in \E(U,V)$.  Alors les   conditions suivantes sont équivalentes :

{\em(i)} pour tout $x\in U$, on a $f(x)\in V^{G^x}$ ;

{\em(ii)} pour tout $x\in U\cap \g r$, on a $f(x)\in V^{G^x}$ ;

{\em(iii)} il existe des applications $f_1,\ldots,f_r\in \E(U)$,
uniquement déterminées, telles
que $f=f_1P_1+\cdots+f_rP_r$.}
\end{prop}
\begin{dem} Il s'agit de montrer (ii)$\Rightarrow$(iii). Soit $x\in \g r$. Il est clair que $V^{\g g^x}=V^{\g g_{\C}^x}$ et 
que $V^{G_{\C}^x}\subset V^{G^x}\subset V^{\g g^x}$, de plus,
comme on l'a déjà signalé, $V^{G_{\C}^x}=V^{\g g_{\C}^x}$, donc
$V^{G^x}=V^{\g g^x}$. Il en découle que pour tout ouvert $U$ de $\g
g$, $f\in \E(U,V)$ (resp. $f\in \A(U,V)$) vérifie la propriété (ii)
(resp. la propriété (ii) de la proposition
\ref{ana}) si et seulement si
\begin{equation}\label{deux}
 \pi(q_i(x))\cdot f(x)=0 \quad\mbox{\rm pour tout  }
i=1,\ldots,\ell \quad\mbox{\rm et pour tout   }x\in U\cap \g r.
\end{equation}
On peut évidemment remplacer dans (\ref{deux}) l'expression $x\in
U\cap \g r$ par l'expression $x\in U$. Considérons les applications
$$ \E(U)^r\stackrel{\iota_U}{\longrightarrow} \E(U,V) \quad\mbox{ et
  }\quad \E(U,V)\stackrel{\tau_U}{\longrightarrow} \E(U,V)^{\ell} $$
définies par  $\iota_U(f_1,\ldots,f_{\ell})=f_1P_1+\ldots +f_rP_r$ et
$\tau_U(f)=(g_1,\ldots,g_{\ell})$, avec $g_i(x)=\pi(q_i(x))\cdot
f(x)$. Tout revient donc à montrer que la suite suivante est exacte.
\begin{equation}\label{trois}
0\longrightarrow \E(U)^r\stackrel{\iota_U}{\longrightarrow}
\E(U,V)\stackrel{\tau_U}{\longrightarrow} \E(U,V)^{\ell}.
\end{equation}

 Notons $\E$
  (resp. $\A$) le faisceau  des germes de fonctions  lisses
  (resp. analytiques) sur $\g g$ à valeurs dans $\C$. Ces notations sont
  cohérentes avec les précédentes, dans le sens où l'ensemble des sections de
  $\E$ au dessus d'un ouvert $U$, qu'on note $\E(U)$,
  est bien l'ensemble des fonctions lisses sur $U$ ;
  on a aussi $(\E\otimes V)(U)=\E(U,V)$. La m\^eme remarque vaut pour
  $\A(U)$ et $\A(U,V)$.

D'après la proposition \ref{ana}, on a l'analogue de (\ref{trois}) pour
tout ouvert $U$ ; d'où la suite exacte de faisceaux
\begin{equation}\label{suiteana}
0\longrightarrow \A^r\stackrel{\iota}{\longrightarrow} \A\otimes
V\stackrel{\tau}{\longrightarrow} (\A\otimes V)^{\ell},
\end{equation}
et, comme $\E$ est fidèlement plat sur $\A$ (\cite{Ma}, chapitre
  VI, corollaire 1.12), en tensorisant  par $\E$, on obtient la suite exacte
$$0\longrightarrow \E^r\stackrel{\iota}{\longrightarrow} \E\otimes
V\stackrel{\tau}{\longrightarrow} (\E\otimes V)^{\ell}.$$
Par un argument bien connu de partition de l'unité, on en déduit
que la suite (\ref{trois}) est exacte.
\end{dem}

\begin{cor}\label{corlisse}{Soient $U$ un ouvert $G$-invariant de $\g g$ et
  $f\in \E(U,V)^{G}$. Alors il existe $r$ fonctions $f_1,\ldots,f_r\in
  \E(U)^{G}$, uniquement déterminées, telles que $f=f_1P_1+\cdots f_rP_r$.}
\end{cor}

\begin{remarque}{On suppose que $V$ admet une forme réelle $V_{\R}$ stable par $G$  et  on choisit les $P_i$ de sorte que leur  restriction à $\g g$
  soit à valeurs dans $V_{\R}$. Alors dans les énoncés \ref{ana},
  \ref{corana}, \ref{lisse} et \ref{corlisse}, si $f$ est à valeurs
  dans $V_{\R}$, les $f_i$ sont à valeurs réelles. 
}
\end{remarque}

\section{Une généralisation d'un théorème de Dixmier}

Dans (\cite{Dix}, théorème 2.6), Dixmier  montre, entre-autres, que, pour tout champ
de vecteurs analytique  $X$ sur $\g g$ tangent aux orbites de $G$,
il existe un champ de vecteurs analytique  $Y$ sur $\g g$ tel que
$X(x)=[x,Y(x)]$ pour tout $x\in\g g$. On se propose d'étendre ce
résultat   aux champs
de vecteurs lisses (c'est une question posée par  Miranda-Zung) et aux
champs de vecteurs à décroissance rapide.

\subsection{Extension aux champs de vecteurs lisses}
On conserve
les notations de la section \ref{secun}. Si $U$ est un ouvert de $\g g$, on note  $\A(U,\g
g)$ l'espace des fonctions  analytiques
sur $U$ à valeurs dans $\g g$.

\begin{prop}\label{lisse}{Soient $U$ un ouvert de $\g g$ et $X\in
    \E(U,\g g)$. Alors les conditions suivantes sont
    équivalentes :

{\em(i)} X annule les fonctions lisses localement invariantes sur $U$
;

{\em(ii)} pour tout $x\in \g r\cap U$, on a $X(x)\in [x,\g g]$ ;

{\em(iii)}  pour tout $x\in U$, on a $X(x)\in [x,\g g]$ ;

{\em(iv)} il existe  une application $Y\in \E(U,\g g)$ 
telle que $X(x)=[x,Y(x)]$ pour tout $x\in U$.
}
\end{prop}
\begin{dem}
Comme dans \cite{Dix}, les implications
(iv)$\Rightarrow$(iii)$\Rightarrow$(ii)$\Rightarrow$(i) et
(i)$\Rightarrow$(ii) sont évidentes, il suffit donc de montrer que
(ii)$\Rightarrow$(iv). La preuve est analogue à celle de la
proposition \ref{lisse}.  On
note  $\A$ (resp. $\E$)  le  faisceau des germes
de fonctions  analytiques  (resp. lisses) sur $\g g$ à valeurs dans $\R$.
Si $\Omega$ est un ouvert de $\g g$ et si $Y\in\A(\Omega,\g g)$,  on définit
$\tau(Y)\in \A(\Omega,\g g)$ par $\tau(Y)(x)=[x,Y(x)]$, et
$\delta(Y)\in \A(\Omega,\g g)^{\ell}$ par
$$\delta(Y)(x)=(\kappa(q_1(x),Y(x)),\ldots,\kappa(q_{\ell}(x),Y(x))),$$
de sorte que, pour $x\in\Omega\cap \g r$, on a $Y(x)\in [x,\g g]$ si
et seulement si $\delta(Y)(x)=0$, car $[x,\g g]$ est l'orthogonal de
$\g g^x$ pour la forme $\kappa$ et
$q_1(x),\ldots,q_{\ell}(x)$ est une base de $\g g^x$.

Le théorème 2.6 de \cite{Dix} affirme que, pour tout ouvert $\Omega$
de $\g g$, la suite 
$$\A(\Omega,\g g)\stackrel{\tau}{\longrightarrow} \A(\Omega,\g
g)\stackrel{\delta}{\longrightarrow}\A(\Omega)^{\ell}$$
est exacte, d'où une  suite exacte de faisceaux
$$\A\otimes\g g \stackrel{\tau}{\longrightarrow}\A\otimes\g g
\stackrel{\delta}{\longrightarrow}\A^{\ell},$$
et  en tonsorisant par $\E$, qui est fidèlement plat sur $\A$, on obtient la suite exacte
$$\E\otimes\g g \stackrel{\tau}{\longrightarrow}\E\otimes\g g
\stackrel{\delta}{\longrightarrow}\E^{\ell} ;$$ 
le résultat en découle par une partition de l'unité.
\end{dem}

En multipliant par des fonctions plateux, on obtient facilement le
corollaire suivant.

\begin{cor}\label{cor1}{Avec les notations de la proposition précédente, si $f$
    est à support compact, on peut choisir $g$ à support compact. }
\end{cor}

\subsection{Extension aux champs de vecteurs à décroissance rapide}

On note $\Sc(\g g,\g g)$ l'espace des fonctions $f\in \E(\g g,\g g)$
vérifiant 
$$\sup_{x\in\g g}\vert Df(x)\vert<+\infty$$
pour tout opérateur différentiel à coefficients polynomiaux $D$ sur
$\g g$.

\begin{prop}\label{schwartz}{Soit $X\in \Sc(\g g,\g g)$.  Alors les   conditions suivantes sont équivalentes :

{\em(i)} X annule les fonctions lisses localement invariantes sur $\g g$
;

{\em(ii)} pour tout $x\in \g r$, on a $X(x)\in [x,\g g]$ ;

{\em(iii)} pour tout $x\in \g g$, on a $X(x)\in [x,\g g]$ ;

{\em(iv)} il existe $Y\in \Sc(\g g,\g g)$ tel que $X(x)=[x,Y(x)]$ pour
tout $x\in\g g$.}
\end{prop}

La démonstration de cette proposition repose sur un lemme de
factorisation inspiré du lemme 8.1.1 de \cite{B}. On fixe une norme
euclidienne $\Vert\cdot \Vert$ sur $\g g$ et un réel $0<a<1$. On note
$B$ la boule unité fermée de $\g g$, $K=\{x\in\g g;
a\leq \Vert x\Vert^2\leq 1\}$ et  $\E_K(\g g,\g g)_0$
l'espace des fonctions $f\in \E_K(\g g,\g g)$  vérifiant
$$f(x)\in [x,\g g]\quad\mbox{\rm pour tout  }x\in\g r.$$

\begin{lemme}\label{factor}{Soit $f\in \E_K(\g g,\g g)_0$.  Il existe une
    nombre fini de fonctions $f_1,\ldots,f_k\in \E_K(\g g,\g
    g)_0$ et des fonctions $\chi_1,\ldots,\chi_r\in \E_{[a,1]}(\R,\R)$
    telles que $f(x)=\chi_1(\Vert x\Vert^2) f_1(x)+\cdots+\chi_k(\Vert
    x\Vert^2) f_k(x)$ pour tout $x\in\g g$.}
\end{lemme}
\begin{dem} La preuve étant  semblable à  celle de (\cite{B}, lemme
    8.1.1), nous nous contentons d'indiquer le point où il y a une
          légère différence.

On note $S$ la sphère unité de $\g g$, $\E(S,\g
g)$ l'espace des fonctions $C^{\infty}$ de $S$ dans $\g g$ et $\E(S,\g
    g)_0$ le sous-espace des $f\in \E(S,\g g)$   vérifiant 
$$f(x)\in [x,\g g]\quad\mbox{\rm pour tout  }x\in S\cap \g r.$$
Il est clair que  $\E(S,\g  g)_0$ est un sous-espace fermé de l'espace de Fréchet $\E(S,\g
g)$, c'est donc un espace de Fréchet pour la topologie induite. Il en
est de m\^eme de l'espace $\E_K(\g g,\g g)_0$. Comme dans la preuve du
lemme 8.1.1 de \cite{B}, on voit que l'application
$$P:\E_K(\g g,\g g)\longrightarrow \E_{[a,1]}(\R , \E(S,\g g))$$
définie par $P(\varphi)(t)(x)=\varphi(t^{\frac 12}x)$ si $t\in [a,1]$,
et $P(\varphi)(t)(x)=0$ si $t\not\in [a,1]$, est un isomorphisme
topologique. Et puisque $\g r$ est un c\^one,   $P$ induit un isomorphisme
topologique de
$\E_K(\g g,\g g)_0$ sur $\E_{[a,1]}(\R ,\E(S,\g
g)_0)$. Le reste de la preuve est
en tout point identique à celui de (\cite{B}, lemme  8.1.1), à ceci
près,  $\E(S,\g g)_0$ remplace  $\E(S)$.
\end{dem}

Revenons à la démonstration de la proposition \ref{schwartz}. On voit,
comme dans le preuve de la proposition \ref{lisse}, qu'il suffit 
de montrer que (ii)$\Rightarrow$(iv). Soit $X\in \Sc(\g g,\g g)$
tel que $X(x)\in [x,\g g]$ pour tout $x\in\g r$. D'après (\cite{Sa},
lemme 2.3), la fonction 
$$\tilde{X}(x)=\left\{\begin{array}{ll}X(\frac{x}{(1-\Vert
  x\Vert^2)^{\frac 12}})& \mbox{\rm si }\Vert x\Vert<1\\0&\mbox{\rm
  sinon }\end{array}\right.$$
appartient à $\E_B(\g g,\g g)$. Il est clair que
  $\tilde{X}(x)\in [x,\g g]$ pour tout $x\in\g r$. 

Soit $\chi\in \E_c(\g g)$ à support dans $B$ et égale à $1$
sur la boule de centre $0$ et de rayon $a$. D'après le corollaire
\ref{cor1}, il existe une fonction $\phi\in \E_c(\g g,\g g)$
telle que $[x,\phi(x)]=\tilde{X}(x)$ pour tout $x\in\g g$ ; d'où
$[x,\chi(x)\phi(x)]=\chi(x)\tilde{X}(x)$ pour tout $x\in\g g$.

La fonction $(1-\chi)\tilde{X}$ est à support dans $K$ et vérifie la
  propriété (ii). Donc, d'après le
  lemme \ref{factor}, on peut écrire
$$(1-\chi)\tilde{X}=\sum_{1\leq i\leq k}(\chi_i\rond \Vert
  \cdot\Vert^2)\cdot f_i,$$
avec $f_i\in \E_K(\g g,\g g)_0$. D'après le corollaire
\ref{cor1}, pour chaque $i$, il existe une fonction $g_i\in
  \E_c(\g g,\g g)$ telle que $[x,g_i(x)]=f_i(x)$ pour tout
  $x\in\g g$. On pose  $\psi=\sum_{1\leq i\leq k}(\chi_i\rond \Vert
  \cdot\Vert^2)\cdot g_i$. Alors  $\psi\in \E_K(\g g,\g g)$
  et $(1-\chi(x))\tilde{X}(x)=[x,\psi(x)]$  pour tout $x\in\g g$.  

La fonction $h=\chi\phi+\psi$ appartient à $ C^{\infty}_B(\g g,\g g)$ et
vérifie $\tilde{X}(x)=[x,h(x)]$  pour tout $x\in\g g$. 

Pour tout $x\in \g g$, on pose 
$$Y(x)=\frac 1{(1+\Vert x\Vert^2)^{\frac 12}}h(\frac x{(1+\Vert
  x\Vert^2)^{\frac 12}}).$$
Il découle alors de la démonstration de (\cite{Sa},
lemme 2.3) que $Y\in \Sc(\g g,\g g)$. On a $X(x)=\tilde{X}(\frac
  x{(1+\Vert x\Vert^2)^{\frac 12}})$ pour tout $x\in \g g$, il découle
  alors de la définition de $Y$ que $X(x)=[x,Y(x)]$ pour tout $x\in\g g$.
{\unskip\hfill\null\nobreak\hfill\carre\vskip1em\par}

\section{Distributions covariantes}

On note $\pi^*$ la représentation duale de $\pi$ dans $V^*$. Chaque
fois qu'on a deux espaces vectoriels en dualité, on note $\langle
\cdot,\cdot\rangle$ le crochet de dualité.

\begin{defi} {Une distribution sur $\g g$ à valeurs dans $V$ est une
    application linéaire continue de $\E_c(\g g)$ dans
    $V$. On note $\D'(\g g,V)$ l'espace de ces
    distributions.}
\end{defi}

L'espace $\D'(\g g,V)$ s'identifie naturellement au dual
topologique de $\E_c(\g g,V^*)$ : si $T\in \D'(\g
g,V)$, on lui associe $S\in \E_c(\g g,V^*)'$, dual topologique
de $\E_c(\g g,V^*)$, définie par
$$\langle S,f\otimes \lambda\rangle =\langle T(f),\lambda\rangle, \quad f\in
\E_c(\g g), \lambda\in V^*.$$

Soient $\chi\in \E(\g g,V)$ et $\theta \in \D'(\g g)$. Si $v_1,\ldots,v_n$ est une base de $V$, on note
$\chi_1,\ldots,\chi_n$ les coefficients de $\chi$ dans cette base, de
sorte que 
$$\chi(x)=\chi_1(x)v_1+\cdots+\chi_n(x)v_n,\quad \mbox{\rm pour
  tout  }x\in\g g.$$
On voit alors facilement qu'on définit une distribution sur $\g g$
à valeurs dans $V$, qu'on notera $\theta \chi$, par
$$\langle \theta  \chi, f\rangle=\sum_{k=1}^n\langle
\theta,f\chi_k\rangle v_k.$$
Cette distribution ne dépend pas du choix  de la base de $V$.

Le groupe $G$ opère dans $\D'(\g g,V)$ par
$$\langle g\cdot T,f\rangle =\pi(g)\cdot \langle T,g^{-1}\cdot f\rangle, \quad g\in G, T\in
\mathcal D'(\g g,V), f\in \E_c(\g g);$$
où $g^{-1}\cdot f$ désigne la fonction $x\mapsto f(\Ad(g^{-1})\cdot
x)$. Une distribution invariante par cette action est appelée
distribution {\it covariante}. On notera $\D'(\g g,V)^G$ le sous-espace de $\D'(\g g,V)$ formé  des distributions  covariantes.

Par exemple si  $\chi\in \E(\g g,V)^G$ et si $\theta\in
\D'(\g g)^G$, la distribution $\theta 
\chi$ appartient à $\D'(\g g,V)^G$.

\begin{remarque}{On rappelle la base $v_1,\ldots,v_n$ de $V$ et on note
    $\eta_1,\ldots,\eta_n$ les fonctions constantes sur $\g g$ à
    valeurs respectives $v_1,\ldots,v_n$. Toute distribution $T\in\D'(\g
    g,V)$ se décompose de façon unique en somme
    $T=T_1\eta_1+\cdots+T_n\eta_n$, avec $T_i\in\D'(\g g)$. Si
    $T$ est covariante, les distributions $T_i$ sont $G$-finies, et,
    réciproquement, toute distribution $G$-finie sur $\g g$ appara\^it
    comme coefficient d'une distribution covariante.}
\end{remarque}

\noindent{\bf Problème 1.} Soit $(P_1,\ldots,P_k)$ une base de
$\P(\g g,V)^G$ comme $\P(\g g)^G$-module. Est-ce que
tout élément de $\D'(\g g,V)^G$ se factorise sous la forme
$$\theta_1 P_1+\cdots+\theta_k P_k, \quad\theta_i\in
\D'(\g g)^G ?$$

Si $V$ admet une forme réelle $V_{\R}$ stable par $G$, il est clair
comment on reformule le problème pour les distributions covariantes à
valeurs dans $V_{\R}$, en choisissant les $P_i$ définis sur $\R$.

En utilisant la méthode de descente de Harish-Chandra, on peut réduire
le problème au cas des distributions covariantes à support inclus dans
le c\^one nilpotent. Mais je n'ai pas réussi à obtenir la
factorisation de ces dernières. La méthode présentée dans la suite
donne, au moins  pour les distributions covariantes à valeurs dans $\g g$, une
décomposition plus faible, et on verra qu'elle fournit des résultats
complets pour $sl(2,\R)$.

 On rappelle la base $(q_1,\ldots,q_{\ell})$ de $\P(\g g,\g
g)^G$ formée de polynômes homogènes (voir paragraphe \ref{troisdeux}), et on identifie $\g g^*$ à $\g g$
 à l'aide de $\kappa$. Dans la suite $\E_c(\g g)$ désigne l'espace des
fonctions réelles de classe $C^{\infty}$ à support compact sur $\g g$,
et $\D'(\g g)$ désigne son dual topologique.

\begin{theo}\label{distcov}{Soit $T\in \D'(\g g,\g g)^G$. Il existe alors
    $\theta_1,\ldots,\theta_{\ell}\in\D'(\g g)$ telles que
    $T=\theta_1q_1+\cdots+\theta_{\ell} q_{\ell}$.}
\end{theo}
\begin{dem} Pour $f\in \E_c(\g g,\g  g)$ et pour 
    $k=1,\ldots,\ell$, on note $\langle f,q_k\rangle$ la fonction
    $x\mapsto \kappa(f(x),q_k(x))$ ; elle appartient à
    $\E_c(\g g)$.

On considère l'application $\psi : \E_c(\g g,\g
    g)\longrightarrow \E_c(\g g)^{\ell}$ définie par
$$\psi(f)=(\langle f,q_1\rangle,\ldots,\langle f,q_{\ell}\rangle).$$
Elle est linéaire continue et l'image
de sa transposée, $^{t}\!\psi$, co\"incide avec l'orthogonal, dans
$\D'(\g g,\g   g)$, du noyau de $\psi$. En effet, par une partition de
l'unité on se ramène au cas des distributions à support compact ; dans
cette situation, le résultat découle  du fait que l'analogue de
l'application $\psi$
 de $\E(\g g,\g g)$ dans $ \E(\g g)^{\ell}$ est à
image fermée (\cite{T}, Chapitre VI,
corollaire 1.5). On pourra trouver une démonstration détaillée dans
(\cite{Sa}, lemme 3.2).

On identifie le dual de  $\E_c(\g g)^{\ell}$ à $\D'(\g
g)^{\ell}$, alors
$^{t}\!\psi(\theta_1,\ldots,\theta_{\ell})=\theta_1
q_1+\cdots+\theta_{\ell} q_{\ell}$. Il suffit donc de montrer
que $T$ appartient à l'orthogonal  du noyau de $\psi$.

Soit $f\in \ker(\psi)$. Pour tout $x\in\g r$, les vecteurs
$q_1(x),\ldots,q_{\ell}(x)$ forment une base de $\g
g^{G^x}=\g g^x$ et l'orthogoanl pour $\kappa$ de $\g g^x$ dans $\g g$ est $[x,\g g]$. Donc, pour tout
$x\in\g r$, $f(x)\in [x,\g g]$. Il découle alors du corollaire
\ref{cor1} qu'il existe une  fonction $\varphi\in
\E_c(\g g)$ telle que $f(x)=[x,\varphi(x)]$ pour tout $x\in \g
g$. 

Si $\xi$ est un élément de $\g g$, on note $\tau(\xi)$ le champ de
vecteurs adjoint sur $\g g$ qui lui est attaché, c'est le champ de
vecteurs défini par $\tau(\xi)(x)=[x,\xi]$ pour tout $x\in\g g$. Il découle alors de la
covariance de $T$ que
\begin{equation}\label{equa5}
\langle T,\ad \xi\rond h\rangle=\langle T,\tau(\xi)\cdot
h\rangle,\quad \mbox{\rm pour tout  }h\in \E_c(\g g, \g g).
\end{equation}

Soit $\xi_1,\ldots,\xi_n$ une base de $\g g$. On note
$\alpha_1,\ldots,\alpha_n$ sa base duale. Alors la
fonction $f$ s'écrit $f=\sum_{i=1}^n\ad \xi_i\rond (\alpha_i
\varphi)$. Donc, d'après (\ref{equa5}), 
\begin{eqnarray*}
\langle T,f\rangle&=&\sum_{i=1}^n\langle T,\ad \xi_i\rond (\alpha_i \varphi)\rangle\\
&=& \sum_{i=1}^n\langle T,\tau(\xi_i)\cdot (\alpha_i \varphi)\rangle\\
&=&\langle T, \sum_{i=1}^n(\tau(\xi_i)\cdot \alpha_i )
\varphi\rangle+\langle T, \sum_{i=1}^n\alpha_i\tau(\xi_i)\cdot \varphi\rangle.
\end{eqnarray*}

On a $\sum_{i=1}^n(\tau(\xi_i)\cdot \alpha_i )(y)=\tr \ad y$ pour tout
$y\in \g g$, donc cette fonction est nulle car $\g g$ est
réductive, et le champ de vecteurs $\sum_{i=1}^n\alpha_i\tau(\xi_i)$
est nul en chaque point, donc $\sum_{i=1}^n\alpha_i\tau(\xi_i)\cdot
\varphi=0$. D'où $\langle T,f\rangle=0$. \end{dem}

On définit de façon évidente l'espace $\mathcal S'(\g g,\g g)^G$ des
distributions tempérées covariantes. Alors, en utilisant la proposition
\ref{schwartz}, on obtient comme ci-dessus :

\begin{theo}\label{distcovt}{Soit $T\in \mathcal  S'(\g g,\g g)^G$. Il existe alors
    $\theta_1,\ldots,\theta_{\ell}\in\mathcal S'(\g g)$ telles que
    $T=\theta_1q_1+\cdots+\theta_{\ell}q_{\ell}$.}
\end{theo}

\begin{remarques}{ 1. Les distributions $\theta_i$ ne sont pas uniques
    et on voudrait  les choisir
    invriantes. Il est clair que chacune d'elles est invariante sur l'image (fermée) de
    l'application $f\mapsto \langle f,q_i\rangle$, et qu'elle n'est
    uniquement déterminée que sur cette image. Pour pouvoir les
    choisir invariantes, il faudrait montrer que toute forme linéaire
    continue invariante sur l'image de $\psi$ se prolonge en une forme
    linéaire continue invariante sur  $\E_c(\g g)^{\ell}$.

2. On pourrait déduire le théorème \ref{distcov} du théorème de division de
   Malgrange (\cite{Ma2},  théorème 1), mais cela
   ne change pas fondamentalement la preuve, puisque la vérification
   des hypothèses du théorème de Malgrange nécessite la description de
   Dixmier   des champs de vecteurs analytiques tangents aux orbites
   et le calcul fait dans la preuve du théorème \ref{distcov}. 

3. Le théorème de Dixmier (\cite{Dix}, théorème 2.6) joue un
    rôle crucial dans la preuve des deux propositions
    précédentes. Pour les étendre aux distributions
    covariantes à valeurs dans n'importe quelle représentation, il
    suffit de résoudre le problème suivant.}
\end{remarques}

\noindent{\bf Problème 2.} On reprend les notations du paragraphe
    \ref{notations} et on rappelle les générateurs
    $q_1,\ldots,q_{\ell}$ de $\mathcal P(\g
    g_{\C},\g g_{\C})^{G_{\C}}$ introduits au début du paragraphe \ref{troisdeux}. Soit $U$ un ouvert de $\g g_{\C}$ qui soit une variété
    de Stein, et soit $f$ une fonction holomorphe de $U$ dans $V$ telle que
    $f(x)$ appartient à l'image de $\pi(x)$ pour tout $x\in \g r\cap
    U$. Existe-t-il des fonctions $g_1,\ldots,g_{\ell}$ holomorphes de $U$
    dans $V$ telles que
$$f(x)=\pi(q_1(x))\cdot g_1(x)+\cdots+\pi(q_{\ell}(x))\cdot g_{\ell}(x),\quad
    \mbox{\rm pour tout }x\in U ?$$

Ce problème a un analogue algébrique, qui pourrait servir pour
appliquer le théorème  de division de Malgrange (voir \cite{Ma3},
paragraphe 3).

\medskip
\noindent{\bf Problème 3.} Soit $P\in\P(\g g_{\C}, V)$  tel que  $P(x)$ appartient à l'image de $\pi(x)$ pour tout $x\in \g  r$. Existe-t-il des polyn\^omes $P_1,\ldots,P_{\ell}\in \P(\g g_{\C}, V)$ tels que
$$P(x)=\pi(q_1(x))\cdot P_1(x)+\cdots+\pi(q_{\ell}(x))\cdot P_{\ell}(x),\quad
    \mbox{\rm pour tout }x\in \g g_{\C}?$$

\section{Cas de $sl(2,\R)$}

On se propose de montrer, dans ce paragraphe,  que la méthode de la
section précédente donne une solution au problème 1 pour $\g g=sl(2,\R)$.

Les représentations irréductibles du groupe adjoint $G$ de $\g g$
ont pour différentielles les représentations irréductibles de
dimension impaire de $\g g$.  Si
$(\pi, V)$ est une représentation irréductible de $G$, l'espace
$\mathcal P(\g g,V)^G$ est un module libre de rang $1$ sur $\mathcal
P(\g g)^G$. On fixe un élément homogène $P$ de $\mathcal P(\g g,V)^G$
tel que $\mathcal P(\g g,V)^G=\mathcal P(\g g)^G P$. Quoique la
représentation contragrédiente de $\pi$ soit équivalente à $\pi$, nous
préférons garder une distinction entre les deux ; on fixe alors un
élément homogène $P^*$ de $\mathcal P(\g g,V^*)^G$ (de m\^eme degré
que $P$) tel que $\mathcal P(\g g,V^*)^G=\mathcal P(\g g)^G
P^*$.

Pour les représentations de dimension impaire de $sl(2,\R)$, le
problème 2 admet une réponse positive. Ici la base
$q_1,\ldots,q_{\ell}$ est formée d'un seul élément, l'identité de $\g g$.

\begin{prop}\label{Dixan1}{ Soit $U$ un ouvert de $\g g_{\C}$ qui soit une variété
    de Stein et soit $f$ une fonction holomorphe de $U$ dans $V^*$ telle
    que $f(x)\in \pi^*(x)(V^*)$ pour tout $x\in \g r_{\C}\cap U$. Alors il
    existe une fonction holomorphe $g$ de $U$ dans $V^*$ telle que
    $f(x)=\pi^*(x)\cdot g(x)$ pour tout $x\in U$.}
\end{prop}
\begin{dem}
On note $E'$ le fibré en droites au dessus de $\g r_{\C}\cap U$ dont
la fibre $E'_x$ en chaque point $x\in \g r_{\C}\cap U$ est égale à
$\ker \pi^*(x)$. Ce fibré admet une section globale, $P^*$, qui ne
s'annule en aucun point de $\g r_{\C}\cap U$, il est donc trivial. La
suite de la preuve est identique à celle du théorème 2.4 de \cite{Dix}.
\end{dem}

\begin{cor}{ Soit $U$ un ouvert de $\g g$ et soit $f\in \A(U,V^*)$
    telle que $f(x)\in \pi^*(x)(V^*)$ pour tout $x\in \g r\cap U$. Alors il
    existe $ \A(U,V^*)$ telle que
    $f(x)=\pi^*(x)\cdot g(x)$ pour tout $x\in U$.}
\end{cor}
\begin{dem} Les arguments de la preuve du théorème 2.6 de \cite{Dix}
  permettent de déduire le corollaire à partir de la proposition
  \ref{Dixan1}.
\end{dem}
\begin{cor}\label{cor2}{ Soit $U$ un ouvert de $\g g$ et soit $f$ un élément de $\E(U,V^*)$
  (resp. $\E_c(U,V^*)$)
    tel que $f(x)\in \pi^*(x)(V^*)$ pour tout $x\in \g r\cap U$. Alors il
    existe $ g\in \E(U,V^*)$ (resp $g\in \E_c(U,V^*)$) tel que
    $f(x)=\pi^*(x)\cdot g(x)$ pour tout $x\in U$.}
\end{cor}
\begin{dem} Résulte du corollaire précédent par les mêmes arguments
    que  la proposition
    \ref{lisse} et le corollaire \ref{cor1}.
\end{dem}

On peut maintenat prouver le résultat principal de ce paragraphe.

\begin{theo}{Soit $T\in \D'(\g g,V)^G$. Alors il existe une distribution
    $\theta\in\D'(\g g)^G$ telle que $T=\theta P$.}
\end{theo}
\begin{dem} On procède comme dans la preuve du théorème
  \ref{distcov}. On considère  l'application $\psi : \E_c(\g
g,V^*)\longrightarrow \E_c(\g g)$ définie par
$\psi(f)(x)=\langle f(x),P(x)\rangle$. L'image de  $\!^t\psi$ co\"incide avec l'orthogonal
du noyau de $\psi$. Le corollaire \ref{cor2} décrit les éléments du
noyau de $\psi$. Un calcul analogue à celui fait dans la preuve de la
proposition  \ref{distcov} montre que $T$ est nulle sur $\ker \psi$,
elle est donc de la forme $SP$ avec $S\in\D'(\g g)$. 

Comme $P(x)\neq
0$ pour tout $x\in \g r=\g g\setminus\{0\}$, l'image de $\psi$
contient $\D(\g g\setminus\{0\})$. On note $E$ l'ensemble des
distributions sur $\g g$ nulles sur l'image de $\psi$ ; c'est un
sous-espace vectoriel de $\D'(\g g)$ stable par $G$ et formé de
distributions à support inclus dans $\{0\}$, donc tout élément de $E$
est $G$-fini.

La restriction de  $S$ à l'image de $\psi$ est $G$-invariante. Donc
pour tout $\xi\in\g g$, la distribtion $\tau(\xi)\cdot S$ est nulle
sur l'image de $\psi$, elle appartient alors à $E$. L'ensemble
des distributions $\tau(\xi)\cdot S$, $\xi\in\g g$, est un sous-espace
vectoriel de dimension finie de $E$. Comme tout
élément de $E$ est $G$-fini, il existe un sous-espace
vectoriel de dimension finie $G$-stable $E_0$ de $E$ tel que $\tau(\xi)\cdot S\in
E_0$ pour tout $\xi\in \g g$.

L'application $\xi\mapsto \tau(\xi)\cdot S$ de $\g g$ dans $E_0$ est clairement un $1$-cocycle. Donc, d'après le lemme de Whitehead, il
existe $\eta\in E_0$ tel que $\tau(\xi)\cdot S=\tau(\xi)\cdot \eta$
pout tout $\xi\in\g g$. La distribution $S-\eta$ est alors $\g
g$-invariante, et donc $G$-invariante car $G$ est connexe. Et puisque
$\eta$ est nulle sur l'image de $\psi$, on a
$\!^t\psi(S-\eta)=\!^t\psi(S)=T$, donc $T=(S-\eta)P$.

\end{dem}

De la m\^eme façon, on obtient l'analogue pour les distributions tempérées.

\begin{theo}{Soit $T\in \mathcal  S'(\g g,V)^G$. Il existe alors
    $\theta\in\mathcal S'(\g g)^G$ telle que
    $T=\theta P$.}
\end{theo}

\begin{remarques}{1. Dans (\cite{lav}, théorème 4.2), Lavaud donne une description des
    distributions covariantes à support inclus dans le c\^one
    nilpotent $\mathcal N$ de $\g g$ légèrement différente de la
    nôtre. Cependant, les arguments concernant l'image
    de $\psi$ dans le preuve du théorème précédent montrent que si le
    support de $T$ est inclus dans $\mathcal N$, alors le support de
    $\theta$ est aussi inclus dans $\mathcal N$. A partir de là, il
    n'est pas difficile de retrouver la formulation de Lavaud. 

2. Si $V$ est de dimension paire, la représentation de $\g g$
   s'intègre  en une représentation de $SL(2,\R)$, mais l'espace des polynômes covariants
de $\g g$ dans $V$ (pour l'action de $SL(2,\R)$) est réduit à $\{0\}$
 ; il en est de m\^eme des différents espaces de fonctions et de
distributions, comme on le voit en regardant l'action du centre de
$SL(2,\R)$. Signalons toutefois que Lavaud a montré que dans ce cas il
existe des distributions localement covariantes.
}
\end{remarques}

\section{Appendice}
Faute de référence, nous présentons ici la preuve de la remarque
\ref{formereelle}. On rappelle les notations : $\g g$  algèbre
de Lie réductive, $\g g_{\C}$ sa complexifiée, $G$ (resp. $G_{\C}$)
groupe adjoint de $\g g$ (resp. $\g g_{\C}$), $(\pi,V)$ une représentation
complexe de dimension finie de $G_{\C}$. 

On suppose que $V$ admet une forme réelle $V_{\R}$ stable par
$G$. Alors $\P(\g g_{\C},V)^{G_{\C}}$ admet une base, comme $\P(\g
g_{\C})^{G_{\C}}$-module, formée d'éléments homogènes $P_1,\ldots,P_r$
tels que $P_i(\g g)\subset V_{\R}$ pour tout $i=1,\ldots,r$.

On note $\sigma_{\g g}$ (resp. $\sigma_{V}$) la conjugaison
de $\g g_{\C}$ (resp. $V$) par rapport  à la forme réelle $\g g$
(resp. $V_{\R}$). Pour $P\in \P(\g g_{\C},V)$, on note $\sigma(P)$
l'élément de $\P(\g g_{\C},V)$ défini par
$$\sigma(P)(x)=\sigma_{V}\cdot P(\sigma_{\g g}\cdot x),\quad\mbox{\rm
  pour tout  }x\in\g g_{\C},$$
et pour $Q\in \P(\g g_{\C})$, on note $\sigma(Q)$
l'élément de $\P(\g g_{\C})$ défini par
$$\sigma(Q)(x)=\overline{Q(\sigma_{\g g}\cdot x)},\quad\mbox{\rm
  pour tout  }x\in\g g_{\C}.$$

Pour  $P\in \P(\g g_{\C},V)$ (resp.  $P\in \P(\g g_{\C})$), on a $P(\g
g)\subset V_{\R}$ (resp. $P(\g g)\subset \R$) si et seulement si
$\sigma(P)=P$ ; on dira dans ce cas que $P$ est défini sur $\R$.

Il s'agit donc de montrer que $\P(\g g_{\C},V)$ admet une base formée
d'élements homogènes définis sur $\R$. On voit par récurrence qu'il suffit de montrer l'assertion suivante.

{\it Pour tout entier $d$, si $\P(\g g_{\C},V)$ admet une base formée
  de polyn\^omes homogènes dont tous les éléments de degré $<d$ sont
  définis sur $\R$, alors il admet une base formée de polyn\^omes homogènes
  dont tous les éléments de degré $\leq d$ sont
  définis sur $\R$.}

Soit  $(Q_1,\ldots,Q_r)$ une base de $\P(\g g_{\C},V)^{G_{\C}}$ formée
d'éléments homogènes et telle que les $Q_i$ de degré $<d$ soient
définis sur $\R$. Quitte à   réindexer les $Q_i$, on peut supposer que
$Q_1,\ldots,Q_k$ sont de degré $d$, $Q_{k+1},\ldots,Q_s$  sont de degré
$<d$ et $Q_{s+1},\ldots,Q_r$  sont  de degré $>d$. Alors pour tout
$j=1,\ldots,k$, il existe des nombres complexes
$\lambda_{1j},\ldots,\lambda_{kj}$ et des polynômes homogènes
$R_{(k+1)j},\ldots,R_{sj}\in P(\g g_{\C})^{G_{\C}}$ tels que
$$\sigma(Q_j)=\sum_{i=1}^k\lambda_{ij}Q_i+\sum_{n=k+1}^sR_{nj}Q_n.$$
D'où
\begin{eqnarray*}
\sigma(\sigma(Q_j))&=&\sigma(\sum_{i=1}^k\lambda_{ij}Q_i+\sum_{n=k+1}^sR_{nj}Q_n)\\
&=&\sum_{i=1}^k\overline{\lambda_{ij}}(\sigma_{m=1}^k\lambda_{mi}Qm+\sum_{n=k+1}^sR_{ni}Q_n)+\sum_{n=k+1}^s\sigma(R_{nj})Q_n\\
&=&\sigma_{m=1}^k(\sum_{i=1}^k\overline{\lambda_{ij}}\lambda_{mi})Q_m+\sum_{n=k+1}^s(\sum_{i=1}^k\overline{\lambda_{ij}}R_{ni})Q_n+\sum_{n=k+1}^s\sigma(R_{nj})Q_n.
\end{eqnarray*}

On note $\Lambda$ la  matrice carrée d'ordre $k$ dont les
coefficients sont les $\lambda_{ij}$, $M$ la matrice à $(s-k)$
lignes et $k$ colonnes dont les coefficients sont les $R_{ij}$ et
$I_k$ la matrice identité de $M_k(\R)$. Il découle alors du calcul
précédent et de $\sigma(\sigma(Q_j))=Q_j$ que
$$\Lambda\overline{\Lambda}=I_k\quad\mbox{et}\quad
R\overline{\Lambda}+\sigma(R)=0.$$
Il existe alors $M=[\mu_{ij}]\in GL(k,\C)$ tel que
$\Lambda=M\overline{M}^{-1}$ (\cite{Se}, chapitre X, proposition
3). D'où
\begin{equation}\label{afin}
\Lambda\overline{M}=M.
\end{equation}
On a
$$\sigma(R\overline{M})=\sigma(R)M=-R\overline{\Lambda}M=-R\overline{M},$$
d'où, en posant $T=\frac 12R\overline{M}$, 
\begin{equation}\label{fin}
 R\overline{M}+\sigma(T)=T.
\end{equation}

On note $T_{ij}$ les coefficients de $T$, et , pour $j=1,\ldots,k$, on pose
$$P_j=\sum_{i=1}^k\mu_{ij}Q_i+\sum_{n=k+1}^sT_{n,j}Q_n.$$
Alors
\begin{eqnarray*}
\sigma(P_j)&=&\sum_{i=1}^k\overline{\mu_{ij}}(\sum_{m=1}^k\lambda_{mi}Q_m+\sum_{n=k+1}^sR_{ni}Q_n)+\sum_{n=k+1}^s\sigma(T_{nj})Q_n\\
&=&\sum_{m=1}^k(\sum_{i=1}^k\overline{\mu_{ij}}\lambda_{mi})Q_m+\sum_{n=k+1}^s[(\sum_{i=1}^k\overline{\mu_{ij}}R_{ni})+\sigma(T_{nj}]Q_n\\
\end{eqnarray*}
Il découle alors de (\ref{afin}) et de (\ref{fin}) que  $\sigma(P_j)=P_j$ pour tout
$1\leq j\leq k$. La famille
$(P_1,\ldots,P_k,Q_{k+1},\ldots,Q_r)$ est donc formée de polyn\^omes homogènes
homogènes et tous ses éléments de degré $\leq d$ sont définis sur
$\R$. Il est
clair qu'elle engendre le $\P(\g
g_{\C})^{G_{\C}}$-module $\P(\g g_{\C},V)^{G_{\C}}$ ; elle est donc une
base (\cite{Bo}, $\S$ 7, n° 10, corollaire 2).

\medskip

\noindent Université de Poitiers-CNRS, Laboratoire de Mathématiques et
Applications,\\ BP 30179, 86962 Futuroscope-Chasseneuil, France.

\noindent {\tt bouaziz@math.univ-poitiers.fr}


\begin{thebibliography}{99}

\bibitem[B]{B}\textsc{A. Bouaziz}  \quad {\em Intégrales orbitales sur
  les algèbres de Lie réductives}, Invent. math. {\bf 115} (1994), 163--207.
\bibitem[Bo]{Bo}\textsc{N. Bourbaki}  \quad {\em Eléménts de
  Mathématiques, Algèbre, Chapitre 2}, Hermann, Paris, 1962.
\bibitem[D]{Dix}\textsc{J. Dixmier}  \quad {\em Champs de vecteurs adjoints
  sur les groupes et algèbres de Lie semi-simples}, J. Reine
  Angew. Math. {\bf 309} (1979), 183--190.
\bibitem[K]{Kos}\textsc{B. Kostant}\quad {\it Lie group
  representations on polynomial rings}, Amer. J. Math. {\bf 85}
  (1963), 327--404.
\bibitem[L]{lav}\textsc{P. Lavaud}\quad {\it Invariant generalized
  functions on $sl(2,\R)$ with values in $sl(2,\R)$}, J. Funct. Anal. {\bf 219}
  (2005), 226--244.
\bibitem[M1]{Ma}\textsc{B. Malgrange}\quad {\it Ideals of
  differentiable functions}, Oxford University Press, 1966.
\bibitem[M2]{Ma2}\textsc{B. Malgrange}\quad {\it Division des
  distributions III : le théorème principal}, Séminaire Schwartz
  (1959-1960),  exposé n° 23-24.
\bibitem[M3]{Ma3}\textsc{B. Malgrange}\quad {\it Division des
  distributions IV : applications}, Séminaire Schwartz
  (1959-1960), exposé n° 25.
\bibitem[MZ]{mir}\textsc{E. Miranda ; N.T. Zung}\quad {\it  A note on
  equivariant normal forms of Poisson structures},  Math. Res. Lett.
  13  (2006), 1001--1012.
\bibitem[O]{ok}\textsc{A.I. Oksak}\quad {\it On invariant and
  covariant Schwartz distributions in the case of a compact linear
  group}, Commun. math. Phys. {\bf 46}
  (1976), 269--287.
\bibitem[S]{Sa}\textsc{A. Saidi}\quad {\it Représentation des
  fonctions et des distributions covariantes sous l'action d'un groupe
  de Lie compact}, Préprint.
\bibitem[Sc]{Schz}\textsc{G.W. Schwarz}\quad {\it Lifting smooth
  homotopies of orbit spaces}, Inst. Hautes Études Sci. Publ. Math. No. 51 (1980), 37--135.
\bibitem[Se]{Se}\textsc{J.-P. Serre}\quad {\it Corps locaux}, Hermann,
  Paris, 1968.
\bibitem[T]{T}\textsc{J. C. Tougeron }\quad {\it Idéaux de fonctions
  différentiables}, Springer-Verlag, Heidelberg, New-York, 1972.
\bibitem[Ve]{Ve}\textsc{F. D. Veldkamp }\quad {\it The center of the
  universal enveloping algebra of a Lie algebra in characteristic
  $p$}, Ann. Scient. E.N.S. {\bf 5}
  (1972), 217--240.



\end{thebibliography}
\end{document}